\documentclass[reqno,12pt]{amsart}

\usepackage{amsmath,amsthm,amssymb,amsfonts,amscd,comment}
\usepackage{xypic}

\theoremstyle{plain} 
\newtheorem{thm}{Theorem}[section]
\newtheorem*{thm*}{Theorem}
\newtheorem{cor}[thm]{Corollary}
\newtheorem{lem}[thm]{Lemma}
\newtheorem{prop}[thm]{Proposition}
\newtheorem*{prop*}{Proposition}

\newtheorem*{conj*}{Conjecture}

\theoremstyle{definition}
\newtheorem{defn}[thm]{Definition}
\newtheorem{eg}[thm]{Example}
\newtheorem{obs}[thm]{Observation}

\theoremstyle{remark}
\newtheorem{rem}[thm]{Remark}

\newtheorem*{pf}{Proof}

\numberwithin{equation}{section}

\def\ZZ{{\mathbb Z}}

\def\RR{{\mathbb R}}
\def\CC{{\mathbb C}}

\def\hh{{\mathbb H}}

\def\XX{{\mathbb X}}

\def\A{{\mathcal A}}
\def\B{{\mathcal B}}
\def\C{{\mathcal C}}
\def\D{{\mathcal D}}

\def\N{{\mathcal N}}

\def\P{{\mathcal P}}

\def\S{{\mathcal S}}
\def\T{{\mathcal T}}

\newcommand{\Aut}{{\rm Aut}}
\newcommand{\Hom}{{\rm Hom}}
\newcommand{\perf}{{\rm perf}}
\newcommand{\Ext}{{\rm Ext}}

\newcommand{\HH}{{\rm H}{\rm H}}

\newcommand{\uSdim}{\overline{{\rm Sdim}}\hspace{0.5mm}}
\newcommand{\lSdim}{\underline{{\rm Sdim}}\hspace{0.5mm}}
\newcommand{\gldim}{{\rm gldim}\hspace{0.5mm}}

\def \mf#1#2#3#4{
\xymatrix{
{#1}\  \ar@<0.4ex>[r]^{{#2}} & \ {#4}
\ar@<0.4ex>[l]^{{#3}}
}
}

\def \mfs#1#2#3#4{\!
\xymatrix@C=1,5em{{#1} \! \ar@<0.2ex>[r]^{{#2}} & \! {#4}
\ar@<0.2ex>[l]^{{#3}}
}
\!}

\def \mfl#1#2#3#4{
\xymatrix@C=2.6em{{#1}\  \ar@<0.4ex>[r]^{{#2}} &\  {#4}
\ar@<0.2ex>[l]^{{#3}}
}
}

\def \mfss#1#2#3#4{\!
\xymatrix@C=1.5em{{#1} \ar@<0.3ex>[r]^{{#2}} & {#4}
\ar@<0.3ex>[l]^{{#3}}
}
\!}
\begin{document}
\title{Serre dimension and stability conditions}
\date{\today}
\author{Kohei Kikuta}
\address{Department of Mathematics, Graduate School of Science, Osaka University, 
Toyonaka Osaka, 560-0043, Japan}
\email{k-kikuta@cr.math.sci.osaka-u.ac.jp}
\author{Genki Ouchi}
\address{Interdisciplinary Theoretical and Mathematical Sciences Program, RIKEN, 2-1 Hirosawa, Wako, Saitama, 351-0198, Japan}
\email{genki.ouchi@riken.jp}
\author{Atsushi Takahashi}
\address{Department of Mathematics, Graduate School of Science, Osaka University, 
Toyonaka Osaka, 560-0043, Japan}
\email{takahashi@math.sci.osaka-u.ac.jp}
\begin{abstract}
We study relations between the Serre dimension defined as the growth of the entropy of the Serre functor
and the global dimension of Bridgeland stability conditions due to Ikeda--Qiu. 
A fundamental inequality between the Serre dimension and the infimum of the global dimensions is proved. 
Moreover, we characterize Gepner type stability conditions on fractional Calabi--Yau categories via the Serre dimension, and classify triangulated categories of Serre dimension lower than one with a Gepner type stability condition. 
\end{abstract}
\maketitle
\section{Introduction}
Dimension is an important notion in mathematics. 
In category theory, 
there has been some attempts to define the dimension of triangulated categories. 
Rouquier defined the dimension (called Rouquier dimension) by the generation-time with respect to a strong generator (\cite{Rou}). 
For autoequivalences of triangulated categories, Dimitrov--Haiden--Katzarkov--Kontsevich defined the notion of entropy motivated by the categorification of classical topological entropy (\cite{DHKK}), which is defined by the growth of generation-time with respect to a split-generator (that is, a classical generator). 
They also computed the entropy of the Serre functor in some cases, and captured a relation to the ``dimension'' of triangulated categories. 
By the computations and comments in \cite{DHKK}, it is natural to come up with a new dimension defined as the growth of the entropy of the Serre functor. 
Recently, Elagin--Lunts defined the upper Serre dimension and the lower Serre dimension in this direction (\cite{EL}).
They also studied basic properties of the Serre dimensions and compared with other notion of dimension (Rouquier dimension, diagonal dimension) of triangulated categories.  
Ikeda--Qiu defined the global dimension of a Bridgeland stability condition on a triangulated category (\cite{IQ,Qiu}), 
which is a natural generalization of the global dimension of finite dimensional algebras. 
They studied in particular the minimal value of global dimensions motivated by the existence of $q$-stability conditions, and observed the equality between the Calabi--Yau dimension and the minimal value of global dimensions in ADE cases. 

In this paper, we introduce the dimensions of triangulated categories defined by the growth of the entropy of the Serre functor, which is exactly the same in \cite{EL}, and use the same terminology i.e. the upper Serre dimension and the lower Serre dimension. 
We also introduce the volume of triangulated categories defined by changing the variable of the entropy of the Serre functor. 
The volume is, by definition, equivalent to the entropy, but this expression gives us a new useful interpretation of the Serre dimension as an analogue of the volume growth of the closed ball in the Euclidean space. 
In other words, the (upper) Serre dimension is considered as a ``similarity dimension" or as a ``scaling dimension". 
It is interesting that the Serre dimension is compatible with the exponent in the theory of Frobenius manifolds or Saito's flat structure, and that Arnold's semicontinuity conjecture is related to the semicontinuity of the Serre dimension of the derived Fukaya category which is homologically mirror to the triangulated category of matrix factorizations of an isolated singularity. 

Moreover, we study relations between the upper Serre dimension $\uSdim\T$ and  the global dimension $\gldim\sigma$ of Bridgeland
stability conditions  $\sigma$.
Firstly, we prove a fundamental inequality.
\begin{thm*}[Theorem \ref{Sdim-le-inf}]
Let $\T$ be a triangulated category equivalent to a perfect derived category of a smooth proper dg $\CC$-algebra. 
We have 
\[
\uSdim\T\le\displaystyle\inf_{\sigma}{\gldim\sigma}. 
\]
Here the infimum runs over all Bridgeland stability conditions on $\T$. 
\end{thm*}
In all known examples, the inequality in the theorem is an equality.
Thus it is natural to ask when $\uSdim\T$ and $\inf_{\sigma}\gldim\sigma$ are equal. 
We give an answer to this question when a fractional Calabi--Yau category admits a Gepner type stability condition.
\begin{thm*}[Theorem \ref{fCY-iff}]
Let $\T$ be a triangulated category equivalent to a perfect derived category of a smooth proper dg $\CC$-algebra. 
Suppose that $\T$ is a fractional Calabi--Yau category with Serre functor $\S$ and a stability condition $\sigma$ on $\T$. 
Then, $\gldim\sigma=\uSdim\T$ if and only if $\sigma$ is of Gepner type with respect to $(\S, \uSdim\T)$ 
(see Definition \ref{def-Gepner}). 
\end{thm*}
The notion of Gepner type stability condition was introduced by Toda, motivated by constructing a stability condition corresponding to the Gepner point of the stringy K\"ahler moduli space of a quintic $3$-fold (\cite{Tod}). 
We note that Gepner type stability is equivalent to $q$-stability. 

Secondly, we study triangulated categories of upper Serre dimension $\le1$. 
Assuming the existence of a stability condition $\sigma$ with $\uSdim\T=\gldim\sigma\le1$, 
the semicontinuity of the upper Serre dimension is proved (Corollary \ref{semiconti}). 
We note that semicontinuity is called monotonicity in \cite{EL}. 
We also classify triangulated categories of $\uSdim<1$ with a Gepner type stability condition.
\begin{thm*}[Theorem \ref{Sdim<1}]
Let $\T$ be a triangulated category equivalent to a perfect derived category of a smooth proper dg $\CC$-algebra. 
Suppose that $\T$ has no nontrivial orthogonal decompositions. 
The following are equivalent.
\begin{enumerate}
\item
$\T$ is equivalent to $\D^b({\rm mod}(\CC Q))$ for some Dynkin quiver $Q$. 
\item
$\uSdim\T<1$ and there exists a Gepner type stability condition on $\T$. 
\item
There exists a stability condition $\sigma$ on $\T$ with $\gldim\sigma<1$. 
\end{enumerate}
\end{thm*}
It would be interesting whether one can remove the condition on the existence of a Gepner type stability condition in (ii). 
As we see in the above statements, the infimum of the global dimensions is also an important invariant. 
Motivated by computations of the infimum in the quiver case due to Qiu, we compute the infimum in  the case of curves.
\begin{thm*}[Theorem \ref{curve-inf}]
Let $C$ be a smooth projective curve of genus $g$.
The following holds. 
\begin{enumerate}
\item
If $g=0$, then there is a stability condition on $\D^b(C)$ such that $\gldim \sigma=1$ and $\displaystyle\inf_{\sigma'}\gldim\sigma'=1$.
\item
If $g=1$, then we have $\gldim \sigma=1$ for any stability condition $\sigma \in {\rm Stab}_\N(\D^b(C))$.
\item
If $g \geq 2$, then $\gldim \sigma>1$ holds for any stability condition  $\sigma \in {\rm Stab}_\N(\D^b(C))$ and $\displaystyle\inf_{\sigma\in {\rm Stab}_\N(\D^b(C))}{\gldim\sigma}=1$.
\end{enumerate}
\end{thm*}
The case of genus greater than one is interesting since there is no minimal value but the infimum is equal to one. 

The contents of this paper are as follows. 
In Section 2, we prepare some notation and define upper and lower Serre dimensions, the Gepner type stability conditions and the global dimension. 
In Section 3, we introduce the volume and $\sigma$-volume of triangulated categories. 
In Section 4, a fundamental inequality (Theorem \ref{Sdim-le-inf}) between the upper Serre dimension and the infimum of the global dimensions is proved. 
Moreover, we prove Theorem \ref{fCY-iff}, which characterizes Gepner type stability conditions on fractional Calabi--Yau categories via the upper Serre dimension.
In Section 5, we study stability conditions with global dimension $\le1$. 
The semicontinuity of the infimum of the global dimensions and of the upper Serre dimension is proved. 
Theorem \ref{Sdim<1} gives the classification for triangulated categories with a Gepner type stability condition whose upper Serre dimension is less than one. 
We also calculate the infimum of global dimensions for derived categories of smooth projective curves (Theorem \ref{curve-inf}). 

\bigskip
\noindent
{\it Acknowledgment}\\
\indent
The authors would like to thank Tom Bridgeland for valuable discussions. 
The first named author is supported by JSPS KAKENHI Grant Number JP17J00227 and the JSPS program ``Overseas Challenge Program for Young Researchers''. 
The second named author is supported by Interdisciplinary Theoretical and Mathematical Science Program (iTHEMS) in RIKEN and JSPS KAKENHI Grant number 19K14520. 
The third named author is supported by JSPS KAKENHI Grant Number 16H06337. 
\section{Preliminaries}
\subsection{Notation}
Throughout this paper, any triangulated category $\T$ is equivalent to $\perf(A)$, 
where $A$ is a smooth proper differential graded (dg) $\CC$-algebra and $\perf(A)$ is the perfect derived category of dg $A$-modules. 
Note that $\T$ has a Serre functor $\S$ and $A$ is a split-generator of $\T$. 
The Grothendieck group (resp. numerical Grothendieck group) of $\T$ is denoted by $K(\T)$ (resp. $\N(\T)$), where $\N(\T)$ is a finitely generated free abelian group by the noncomuutative Hirzebruch--Riemann--Roch theorem (\cite{Shk, Lun}). 
We denote the group of autoequivalences of $\T$ by $\Aut(\T)$. 

A variety means an integral separated scheme of finite type over $\CC$. 
Throughout this paper, $X$ is a smooth projective variety over $\CC$ and $Q$ is a finite connected quiver. 
The bounded derived category of coherent sheaves on $X$ is denoted by $\D^b(X):=\D^b({\rm Coh}(X))$. 
For a finite dimensional $\CC$-algebra $B$ (resp. a quiver $Q$), the bounded derived category of finitely generated $B$-modules (resp. $\CC Q$-modules) is denoted by $\D^b(B):=\D^b({\rm mod}\hspace{0.5mm}(B))$ (resp. $\D^b(Q):=\D^b({\rm mod}\hspace{0.5mm}(\CC Q))$).  
\subsection{Serre dimension}
In this subsection, we give the definition of the Serre dimension of triangulated categories. 
\begin{defn}[{\cite[Definition 2.5]{DHKK}}]
Let $G\in\T$ be a split-generator and $F\in\Aut(\T)$. 
The {\it entropy} of $F$ is the function $h_t(F):\RR\to[-\infty,+\infty)$ defined by
\begin{equation*}
h_{t}(F):=\displaystyle\lim_{n\rightarrow+\infty}\frac{1}{n}\log \delta_t(G,F^{n}G), 
\end{equation*}
where
\begin{equation*}
\delta_t(G,F^{n}G)
:=
\inf\left\{
\displaystyle\sum_{i=1}^p {\rm exp}(n_i t)~
\middle
|~
\begin{xy}
(0,5) *{0}="0", (20,5)*{M_{1}}="1", (30,5)*{\dots}, (40,5)*{M_{p-1}}="k-1", (60,5)*{F^nG\oplus M}="k",
(10,-5)*{G[n_{1}]}="n1", (30,-5)*{\dots}, (50,-5)*{G[n_{p}]}="nk",
\ar "0"; "1"
\ar "1"; "n1"
\ar@{.>} "n1";"0"
\ar "k-1"; "k" 
\ar "k"; "nk"
\ar@{.>} "nk";"k-1"
\end{xy}
\right\}
\end{equation*}
\end{defn}
In the definition of the entropy, the limit exists and doesn't depend on the choice of split-generators (see \cite[Lemma 2.6]{DHKK}). 
\begin{thm}[{\cite[Theorem 2.7]{DHKK}}]\label{DHKK-thm}
Let $G$ be a split-generator of $\T$ and $F\in\Aut(\T)$. 
The entropy $h_t(F)$ is given by
\begin{equation*}
h_t(F)=\lim_{n\rightarrow+\infty}\frac{1}{n}\log\sum_{m\in\ZZ}\dim_{\CC} {\rm Hom}_\T(G,F^nG[m]) e^{-mt}. 
\end{equation*}
\end{thm}
The following lemma is a direct corollary of this theorem. 
\begin{lem}\label{serre-lim-exist}
The limits $\displaystyle\lim_{t\to+\infty}{\frac{h_t(\S)}{t}}$ and $\displaystyle\lim_{t\to-\infty}{\frac{h_t(\S)}{t}}$ exist, especially we have
\begin{equation*}
\displaystyle
\lim_{t\to+\infty}{\frac{h_t(\S)}{t}}=\limsup_{n\to+\infty}{\frac{m^-_n}{n}},~
\displaystyle\lim_{t\to-\infty}{\frac{h_t(\S)}{t}}=\limsup_{n\to+\infty}{\frac{m^+_n}{n}}, 
\end{equation*}
where, for each $n\in\ZZ_{>0}$, set
\begin{eqnarray*}
m^-_n=m^-_n(G):=-\min\left\{m\in\ZZ~|~\Hom_\T(G,\S^n G[m])\neq0\right\}\\
m^+_n=m^+_n(G):=-\max\left\{m\in\ZZ~|~\Hom_\T(G,\S^n G[m])\neq0\right\}. 
\end{eqnarray*}
\end{lem}
\begin{pf}
Recall that $\S$ is the Serre functor of $\T$. We prove the first equality. 
Let $G$ be a split-generator of $\T$. 
It is easy to see that
\begin{eqnarray*}
e^{m^-_nt}\dim_{\CC} {\rm Hom}_\T(G,\S^nG[-m^-_n])
&\le&\sum_{m\in\ZZ}\dim_{\CC} {\rm Hom}_\T(G,\S^nG[m]) e^{-mt}\\
&\le&e^{m^-_nt}\sum_{m\in\ZZ}\dim_{\CC} {\rm Hom}_\T(G,\S^nG[m])
\end{eqnarray*}
for $t\ge1$. By Theorem \ref{DHKK-thm}, we obtain the first equality. 
The second equality follows from the same argument. 
\qed
\end{pf}
\begin{defn}\label{def of Serre dim}
The {\it upper Serre dimension} $\uSdim\T\in[-\infty,+\infty]$ of $\T$ is given by
\begin{equation*}
\uSdim\T:=\lim_{t\to+\infty}{\frac{h_t(\S)}{t}}. 
\end{equation*}
The {\it lower Serre dimension} $\lSdim\T\in[-\infty,+\infty]$ of $\T$ is given by
\begin{equation*}
\lSdim\T:=\lim_{t\to-\infty}{\frac{h_t(\S)}{t}}. 
\end{equation*}
\end{defn}
It is clear that $\uSdim\T\ge\lSdim\T$ by the definition. 
\begin{eg}\label{example-dimension}
The following are examples of the Serre dimensions. 
\begin{enumerate}
\item
$\uSdim\D^b(X)=\lSdim\D^b(X)=\dim X$. 
\item
$\uSdim\D^b(Q)=\lSdim\D^b(Q)=1$, where $Q$ is an acyclic quiver, not of Dynkin type (see \cite[Theorem 2.17]{DHKK}). 
\item
$\uSdim\T=\lSdim\T=\frac{k_2}{k_1}$, where $\T$ is a fractional Calabi--Yau category such that $\S^{k_1}=[k_2]$.
\item 
Let $Q$ be a Dynkin quiver. Then $\D^b(Q)$ is a fractional Calabi-Yau category such that $\S^h=[h-2]$, where $h$ is the Coxeter number of $Q$ (see \cite[8.3 (2)]{Kel}).
\end{enumerate}
\end{eg}
\begin{prop}[{\cite[Proposition 6.14]{EL}}]\label{ob1}
Fix $d\in\RR$. 
Then $\uSdim\T=\lSdim\T=d$ if and only if $h_t(\S)=d\cdot t+h_0(\S)$. 
\end{prop}
\begin{rem}
From the view point of Frobenius structures or Saito's flat structures
(cf. \cite{Dub, ST}),
it seems more important to consider the dimension in complex numbers,
by taking the imaginary part into account in addition to the Serre dimension as the real part. 
Namely, if $h_t(\S)=d\cdot t+h_0(\S)$, then it is natural to introduce 
\begin{equation*}
\hat{c}_\T:=d+\sqrt{-1}~\frac{h_0(\S)}{\pi}=\frac{``h_{-\pi\sqrt{-1}}(\S)"}{-\pi\sqrt{-1}}\in\CC.
\end{equation*}
In the theory of Frobenius structures, the dimension is defined as the difference between the largest and 
the smallest exponents/spectral numbers, those shall be related to logarithms of eigenvalues of 
the automorphism on the Grothendieck group induced by the Serre functor.
For a triangulated category $\T$ admitting a full exceptional collection, there is an idea to reconstruct exponents/spectral numbers 
from the Euler form $\chi$, which was proposed by Cecotti--Vafa \cite{CV} and was developed by Balnojan--Hertling \cite{BH}. 
Let us give a simple example to explain this relation.
If $Q$ is the Kronecker quiver with more than two arrows,  then by \cite{CV, BH} one obtains a complex number
$\hat{c}_\T$ from the Euler form $\chi$ for $\T=\D^b(Q)$. In particular, 
${\rm exp}(\pm\pi\sqrt{-1}\hat{c}_\T)$ are eigenvalues of $\chi^{-1}\chi^T$, the automorphism on the Grothendieck group 
induced by the Serre functor. This complex number $\hat{c}_\T$ is the dimension of the semi-simple Frobenius structure 
of rank two whose Stokes matrix gives rise to the Euler form $\chi$. 
\end{rem}
\subsection{Global dimension of stability conditions}
In this subsection, we give the definition of the global dimension function of Bridgeland stability conditions, due to Ikeda--Qiu. 

Fix a finitely generated free abelian group $\Lambda$, a surjective group homomorphism $v: K(\T)\twoheadrightarrow\Lambda$ and a group homomorpshim $\alpha: \Aut(\T)\to\Aut_{\ZZ}(\Lambda)$, such that the following diagram of abelian groups is commutative:
\[
  \xymatrix{
    K(\T) \ar[r]^{K(F)} \ar[d]_{v} & K(\T) \ar[d]_{v} \\
    \Lambda \ar[r]^{\alpha(F)} & \Lambda 
  }
\]
for any $F \in \Aut(\T)$.
We also fix a norm $||\cdot||$ on $\Lambda\otimes_\ZZ \RR$. 
\begin{defn}[{\cite[Definition 5.1]{Bri}}]\label{stab-def}
A {\it stability condition} $\sigma=(Z,\P)$ on $\T$ (with respect to $(\Lambda, v)$) consists of a group homomorphism
$Z: \Lambda\to\CC$ called {\it central charge} and a family $\P=\{\P(\phi)\}_{\phi\in\RR}$ of full additive subcategories of $\T$ called {\it slicing}, such that
\begin{enumerate}
\item
For $0\neq E\in\P(\phi)$, we have $Z(v(E))=m(E)\exp(i\pi\phi)$ for some $m(E)\in\RR_{>0}$. 
\item
For all $\phi\in\RR$, we have $\P(\phi+1)=\P(\phi)[1]$.
\item
For $\phi_1>\phi_2$ and $E_i\in\P(\phi_i)$, we have $\Hom(E_1,E_2)= 0$.
\item
For each $0\neq E\in\T$, there is a collection of exact triangles called {\it Harder--Narasimhan filtration} of $E$: 
\begin{equation}\label{HN}
\begin{xy}
(0,5) *{0=E_0}="0", (20,5)*{E_{1}}="1", (30,5)*{\dots}, (40,5)*{E_{p-1}}="k-1", (60,5)*{E_p=E}="k",
(10,-5)*{A_{1}}="n1", (30,-5)*{\dots}, (50,-5)*{A_{p}}="nk",
\ar "0"; "1"
\ar "1"; "n1"
\ar@{.>} "n1";"0"
\ar "k-1"; "k" 
\ar "k"; "nk"
\ar@{.>} "nk";"k-1"
\end{xy}
\end{equation}
with $A_i\in\P(\phi_i)$ and $\phi_1>\phi_2>\cdots>\phi_p$. 
\item
(support property)
There exists a constant $C>0$ such that for all $0\neq E\in\P(\phi)$, 
we have 
\begin{equation}\label{supp-prop}
||v(E)||<C|Z(v(E))|.
\end{equation}
\end{enumerate}
\end{defn}
For any interval $I\subset\RR$, define $\P(I)$ to be the extension-closed subcategory of $\T$ generated by the subcategories $\P(\phi)$ for $\phi\in I$. 
Then $\P((0,1])$ is the heart of a bounded t-structure on $\T$, hence an abelian category.  
The full subcategory $\P(\phi)\subset\T$ is also shown to be abelian. 
A nonzero object $E\in\P(\phi)$ is called {\it $\sigma$-semistable} of {\it phase} $\phi_\sigma(E):=\phi$, and especially a simple object in $\P(\phi)$ is called {\it $\sigma$-stable}. 
Taking the Harder--Narasimhan filtration (\ref{HN}) of $E$, we define $\phi^+_\sigma(E):=\phi_\sigma(A_1)$ and $\phi^-_\sigma(E):=\phi_\sigma(A_p)$. 
The object $A_i$ is called {\it $\sigma$-semistable factor} of $E$. 
Define ${\rm Stab}_\Lambda(\T)$ to be the set of stability conditions on $\T$ with respect to $(\Lambda,v)$, 
especially ${\rm Stab}_\N(\T)$ to be the set of stability conditions on $\T$ with respect to the natural projection $K(\T)\to\N(\T)$. 
An element in ${\rm Stab}_\N(\T)$ is called a {\it numerical stability condition} on $\T$. 

In this paper, we assume that the space ${\rm Stab}_\Lambda(\T)$ is not empty for some $(\Lambda,v)$. 
We will abuse notation and write $Z(E)$ instead of $Z(v(E))$. 

We prepare some terminologies on the stability on the heart of a $t$-structure on $\T$. 
\begin{defn}
Let $\A$ be the heart of a bounded $t$-structure on $\T$. 
A {\it stability function} on $\A$ is a group homomorphism $Z: \Lambda\to\CC$ such that for all $0\neq E\in\A\subset\T$, the complex number $Z(v(E))$ lies in the semiclosed upper half plane $\hh_{-}:=\{re^{i\pi\phi}\in\CC~|~r\in\RR_{>0},\phi\in(0,1]\}\subset\CC$. 
\end{defn}
Given a stability function $Z: \Lambda\to\CC$ on $\A$, the {\it phase} of an object $0\neq E\in\A$ is defined to be $\phi(E):=\frac{1}{\pi}{\rm arg}Z(E)\in(0,1]$. 
An object $0\neq E\in\A$ is {\it $Z$-semistable} (resp. {\it $Z$-stable}) if for all subobjects $0\neq A\subset E$, we have $\phi(A)\le\phi(E)$ (resp. $\phi(A)<\phi(E)$). 
We say that a stability function $Z$ satisfies {\it the Harder--Narasimhan property} if 
each object $0\neq E\in\A$ admits a filtration (called Harder--Narasimhan filtration of $E$) 
$0=E_0\subset E_1\subset E_2\subset\cdots\subset E_m=E$ such that $E_i/E_{i-1}$ is $Z$-semistable for $i=1,\cdots,m$ with $\phi(E_1/E_0)>\phi(E_2/E_1)>\cdots>\phi(E_m/E_{m-1})$.
A stability function $Z$ on $\A$ satisfies {\it the support property} if there exists a constant $C>0$ such that for all $Z$-semistable objects $E\in\A$, 
we have $||v(E)||<C|Z(v(E))|$. 

The following proposition shows the relationship between stability conditions and stability functions on the heart of a bounded $t$-structure. 
\begin{prop}[{\cite[Proposition 5.3]{Bri}}]
To give a stability condition on $\T$ is equivalent to giving a bounded t-structure on $\T$ with the heart $\A$, and a stability function $Z$ on $\A$ with the Harder--Narasimhan property and the support property. 
\end{prop}
For the proof, we construct the slicing $\P$, from the pair $(Z,\A)$, by 
\[
\P(\phi):=\{E\in\A~|~E\text{ is }Z\text{-semistable with }\phi(E)=\phi\}\text{ for }\phi\in(0,1],
\]
and extend for all $\phi\in\RR$ by $\P(\phi+1):=\P(\phi)[1]$. 
Conversely, for a stability condition $\sigma=(Z,\P)$, the heart $\A$ is given by $\A:=\P_\sigma((0,1])$. 
We also denote stability conditions by $(Z,\A)$. 

There are two natural group-actions on ${\rm Stab}_\Lambda(\T)$. 
The first is the left $\Aut(\T)$-action defined by 
\begin{equation*}
F.\sigma:=(Z_\sigma(\alpha(F^{-1})(-)), \{F(\P_\sigma(\phi))\})\text{ for }F\in\Aut(\T). 
\end{equation*}
The second is the right $\CC$-action defined by 
\begin{equation*}
\sigma.\mu:=(\exp(-i\pi\mu)\cdot Z_\sigma, \{\P_\sigma(\phi+{\rm Re}(\mu))\})\text{ for }\mu\in\CC. 
\end{equation*}
The notion of the Gepner type stability condition was introduced by Toda, motivated by constructing a stability condition corresponding to the Gepner point of the stringy K\"ahler moduli space of a quintic $3$-fold. 
This notion plays a central role in Section \ref{sec-stab-Gepner}. 
\begin{defn}[{\cite[Definition 2.3]{Tod}}]\label{def-Gepner}
A stability condition $\sigma$ on $\T$ is {\it Gepner type} with respect to $(F,\mu)\in\Aut(\T)\times\CC$ 
if the condition $F.\sigma=\sigma.\mu$ holds. 
\end{defn}
\begin{eg}[{\cite[Theorem 4.2]{KST1}(see also \cite[Theorem 2.14]{Tod})}]\label{KST-Gepner}
Kajiura--Saito--Takahashi constructed the Gepner type stability condition $\sigma_G$ $\in {\rm Stab}_\N(\D^b(Q))$
with respect to $(\S,1-\frac{2}{h_Q})$, where $Q$ is a Dynkin quiver and $h_Q>0$ is the Coxeter number of $Q$. 
\end{eg}
The global dimension was introduced by Ikeda--Qiu for analyzing $q$-stability conditions.
This notion is a natural generalization of the global dimension of finite dimensional algebras. 
The Serre dimension is closely related to (the infimum of) the global dimension function (see Theorem \ref{Sdim-le-inf}). 
\begin{defn}[{\cite[Definition 2.20]{IQ}}]
For a stability condition $\sigma=(Z,\P)$ on $\T$, 
the {\it global dimension} $\gldim\sigma$ of $\sigma$ is given by
\begin{equation*}
\gldim\sigma:=\sup\left\{\phi_2-\phi_1|~\Hom_\T(E_1,E_2)\neq0\text{ for }E_i\in\P(\phi_i)\right\}\in[0,+\infty]. 
\end{equation*}
\end{defn}
The global dimension function is continuous with respect to some natural topology on ${\rm Stab}_\Lambda(\T)$(\cite[Lemma 2.23]{IQ}).  
\section{Volume}
We introduce the notion of the volume and the $\sigma$-volume of triangulated categories, and study the relation to the Serre dimensions. 
\subsection{Volume}
We introduce the volume of triangulated categories via the entropy of the Serre functor. 
\begin{defn}
For $\lambda\in\RR_{>0}$, the {\it volume} $V_\lambda(\T)$ of $\T$ at scale $\lambda$ is defined by
\begin{equation*}
V_\lambda(\T):=\exp(h_{\log\lambda}(\S))\in[0,+\infty). 
\end{equation*}
\end{defn}
The following is an important observation. 
\begin{obs}
Let $(\RR^d,g_E)$ be the $d$-dimensional Euclidean space and $\bar{B}(\lambda)$ the closed ball in $\RR^d$ with center the origin and radius $\lambda>0$. 
Then the volume ${\rm Vol}_\lambda$ of $\bar{B}(\lambda)$ is $\displaystyle{\rm Vol}_\lambda=\frac{\pi^{\frac{d}{2}}}{\Gamma(\frac{d}{2}+1)}\lambda^d$. 
Therefore we have ${\rm Vol}_\lambda={\rm Vol}_1\cdot\lambda^d$, and the dimension $d$ is described by
\begin{equation*}
d=\displaystyle\lim_{\lambda\to+\infty}\frac{\log {\rm Vol}_\lambda}{\log\lambda}
=\displaystyle\lim_{\lambda\to+0}\frac{\log {\rm Vol}_\lambda}{\log\lambda}.  
\end{equation*}
\end{obs}
The following is clear by the definition of the volume. 
This however gives us an another useful interpretation of the Serre dimensions motivated by the above observation. 
\begin{prop}\label{serre-vol}
We have the following. 
\begin{enumerate}
\item
$\uSdim\T=\displaystyle\lim_{\lambda\to+\infty}\frac{\log V_\lambda(\T)}{\log\lambda}$. 
\item
$\lSdim\T=\displaystyle\lim_{\lambda\to+0}\frac{\log V_\lambda(\T)}{\log\lambda}$. 
\end{enumerate}
\end{prop}
Therefore the (upper) Serre dimension can be interpreted as a ``similarity dimension" or as a ``scaling dimension". 

By the Proposition \ref{ob1}, we obtain the following. 
\begin{cor}\label{ob2}
Fix $d\in\RR$. 
Then $\uSdim\T=\lSdim\T=d$ if and only if  
\begin{equation*}
V_\lambda(\T)=V_1(\T)\cdot\lambda^d. 
\end{equation*}
for any $\lambda>0$.
\end{cor}
Therefore the equality between the upper Serre dimension and the lower Serre dimension seems to be a natural condition (see also Proposition \ref{ob1}). 
\subsection{$\sigma$-Volume}
We introduce the $\sigma$-volume of triangulated categories via the mass-growth with respect to the Serre functor, and show some results similar to that of the volume. 
The purpose of this subsection is to propose the analogue of the volume via the mass-growth, 
thus for the proof of the main theorems, one can skip to the next section. 
\begin{defn}
Let $E\in\T$ be a nonzero object of $\T$ and $\sigma\in {\rm Stab}_\Lambda(\T)$ be a stability condition on $\T$. 
The {\it mass} of $E$ with a parameter $t\in\RR$ is the function $m_{\sigma,t}(E):\RR\to\RR_{>0}$ defined by
\begin{equation*}
m_{\sigma,t}(E):=\displaystyle\sum_{i=1}^p|Z_\sigma(A_i)|e^{\phi_\sigma(A_i)t},
\end{equation*}
where $A_1,\cdots,A_p$ are $\sigma$-semistable factors of $E$. 
\end{defn}
\begin{defn}[{\cite[Section 4]{DHKK} and \cite[Theorem 3.5(1)]{Ike}}]
Let $G\in\T$ be a split-generator, $F\in\Aut(\T)$ an autoequivalence of $\T$ and 
$\sigma\in {\rm Stab}_\Lambda(\T)$ a stability condition on $\T$. 
The {\it mass-growth} with respect to $F$ is the function $h_{\sigma, t}(F):\RR\to[-\infty,+\infty]$ defined by
\begin{equation*}
h_{\sigma,t}(F):=\displaystyle\limsup_{n\rightarrow+\infty}\frac{1}{n}\log m_{\sigma,t}(F^n G)
\end{equation*}
It does not depend on a choice of a split generator $G$.
\end{defn}
\begin{thm}[{\cite[Thorem 3.5(2)]{Ike}}]\label{mass-growth-le-entropy}
Let $F\in\Aut(\T)$ be an autoequivalence of $\T$ and $\sigma\in {\rm Stab}_\Lambda(\T)$ a stability condition on $\T$. 
Then we have
\begin{equation*}
h_{\sigma,t}(F)\le h_t(F)<+\infty.
\end{equation*}
\end{thm}
\begin{defn}
For $\lambda\in\RR_{>0}$ and $\sigma\in {\rm Stab}_\Lambda(\T)$, the {\it $\sigma$-volume} $V_{\sigma,\lambda}(\T)$ of $\T$ at scale $\lambda$ is defined by
\begin{equation*}
V_{\sigma,\lambda}(\T):=\exp(h_{\sigma,\log\lambda}(\S))\in[0,+\infty). 
\end{equation*}
\end{defn}
\begin{prop}\label{lim-of-phase}
Let $G$ be a split-generator of $\T$, $\sigma\in {\rm Stab}_\Lambda(\T)$ a stability condition on $\T$. 
Then we have the following.
\begin{enumerate}
\item
$\displaystyle\uSdim\T
=\limsup_{n\to+\infty}\frac{1}{n}\phi^+_\sigma(\S^n G)
=-\limsup_{n\to+\infty}\frac{1}{n}\phi^-_\sigma(\S^{-n} G)$. 
\item
$\displaystyle\lSdim\T
=\limsup_{n\to+\infty}\frac{1}{n}\phi^-_\sigma(\S^n G)
=-\limsup_{n\to+\infty}\frac{1}{n}\phi^+_\sigma(\S^{-n} G)$. 
\end{enumerate}
\end{prop}
\begin{pf}
We prove (i). 
Since $\Hom(G,\S^nG[-m^-_n])\neq0$ ($n\in\ZZ_{>0}$), we have 
\begin{equation*}
\phi^+_\sigma(\S^nG)-\phi^-_\sigma(G)\ge m^-_n\text{ and }\phi^+_\sigma(G)-\phi^-_\sigma(\S^{-n}G)\ge m^-_n,
\end{equation*}
which imply $\displaystyle\limsup_{n\to+\infty}\frac{1}{n}\phi^+_\sigma(\S^n G)\ge\uSdim\T$ and $\displaystyle-\limsup_{n\to+\infty}\frac{1}{n}\phi^-_\sigma(\S^{-n} G)\ge\uSdim\T$. 
Now we prove the opposite inequalities. 
Let $A^{(n)}$ be the $\sigma$-semistable factor of $\S^nG$ with $\phi_\sigma(A^{(n)})=\phi^+_\sigma(\S^nG)$, 
and $A'^{(n)}$ be the $\sigma$-semistable factor of $\S^{-n}G$ with $\phi_\sigma(A'^{(n)})=\phi^-_\sigma(\S^{-n}G)$. 

We can take $l\in\ZZ_{>0}$ satisfying $\displaystyle l\cdot\frac{C'}{C}\ge1$, 
where $C>0$ is the constant appearing in the support property (\ref{supp-prop}), 
and 
$C':=\min\{||\gamma||~|~\gamma\in\Lambda\backslash\{0\}\}>0$. 
For the split-generator $G_l:=G\oplus\cdots\oplus G$ ($l$-th direct sum), it is worth to note that 
$m_n^-(G_l)=m_n^-(G), 
\phi^+_\sigma(\S^nG_l)=\phi^+_\sigma(\S^nG)$
 and $\phi^-_\sigma(\S^{-n}G_l)=\phi^-_\sigma(\S^{-n}G)$.  
The $l$-th direct sum $A^{(n)}_l:=A^{(n)}\oplus\cdots\oplus A^{(n)}$  (resp. $A'^{(n)}_l:=A'^{(n)}\oplus\cdots\oplus A'^{(n)}$)
is the $\sigma$-semistable factor of $\S^nG_l$ with $\phi_\sigma(A^{(n)}_l)=\phi^+_\sigma(\S^nG_l)$
(resp. the $\sigma$-semistable factor of $\S^{-n}G_l$ with $\phi_\sigma(A'^{(n)}_l)=\phi^-_\sigma(\S^{-n}G_l)$). 
By the support property, we have 
\begin{eqnarray*}
|Z_\sigma(A^{(n)}_l)|=l\cdot|Z_\sigma(A^{(n)})|
\ge l\cdot\frac{C'}{C}\ge1\\
|Z_\sigma(A'^{(n)}_l)|=l\cdot|Z_\sigma(A'^{(n)})|
\ge l\cdot\frac{C'}{C}\ge1
\end{eqnarray*}

Then by the definition of the mass-growth, we have 
\begin{eqnarray*}
h_{\sigma,t}(\S)
&\ge&t\cdot\limsup_{n\to+\infty}\frac{1}{n}\phi^+_\sigma(\S^nG_l)
=t\cdot\limsup_{n\to+\infty}\frac{1}{n}\phi^+_\sigma(\S^nG)\\
h_{\sigma,-t}(\S^{-1})
&\ge&(-t)\cdot\limsup_{n\to+\infty}\frac{1}{n}\phi^-_\sigma(\S^{-n}G_l)
=(-t)\cdot\limsup_{n\to+\infty}\frac{1}{n}\phi^-_\sigma(\S^{-n}G). 
\end{eqnarray*}
Here we use the elementary inequality
$
\limsup_{n\to+\infty}(a_n+b_n)\ge\limsup_{n\to+\infty}b_n
$
for sequences $\{a_n\}_n$ and $\{b_n\}_n$ such that $a_n\ge0$ for all $n$. 
Therefore by Theorem \ref{mass-growth-le-entropy}, we have
\begin{eqnarray*}
\uSdim\T=\lim_{t\to+\infty}\frac{h_t(\S)}{t}\ge\limsup_{t\to+\infty}\frac{h_{\sigma,t}(\S)}{t}\ge\limsup_{n\to+\infty}\frac{1}{n}\phi^+_\sigma(\S^nG)\\
\uSdim\T=\lim_{t\to+\infty}\frac{h_t(\S)}{t}=\lim_{t\to+\infty}\frac{h_{-t}(\S^{-1})}{t}\ge\limsup_{t\to+\infty}\frac{h_{\sigma,-t}(\S^{-1})}{t}\ge-\limsup_{n\to+\infty}\frac{1}{n}\phi^-_\sigma(\S^{-n}G). 
\end{eqnarray*}
The statement (ii) follows from the same argument. 
\qed
\end{pf}
\begin{lem}\label{sigma-vol-scale}
We have the following.
\begin{enumerate}
\item
$V_{\sigma,1}(\T)\cdot\lambda^{\lSdim\T}\le V_{\sigma,\lambda}(\T)\le V_{\sigma,1}(\T)\cdot\lambda^{\uSdim\T}\text{  for }\lambda\ge1$. 
\item
$V_{\sigma,1}(\T)\cdot\lambda^{\uSdim\T}\le V_{\sigma,\lambda}(\T)\le V_{\sigma,1}(\T)\cdot\lambda^{\lSdim\T}
\text{  for }0<\lambda<1$. 
\end{enumerate}
\end{lem}
\begin{pf}
We prove (i). 
By the definition of the mass-growth, it is easy to see that
\begin{equation*}
\lambda^{\phi^-_\sigma(\S^nG)}m_{\sigma,0}(\S^nG)
\le m_{\sigma,\log\lambda}(\S^nG)
\le \lambda^{\phi^+_\sigma(\S^nG)}m_{\sigma,0}(\S^nG),
\end{equation*}
which gives the inequality by Proposition \ref{lim-of-phase}. 
The statement (ii) follows from the same argument. 
\qed
\end{pf}
\begin{lem}\label{serre-sigma-vol}
We have the following.
\begin{enumerate}
\item
$\uSdim\T=\displaystyle\lim_{\lambda\to+\infty}\frac{\log V_{\sigma,\lambda}(\T)}{\log\lambda}=\lim_{t\to+\infty}\frac{h_{\sigma,t}(\S)}{t}$. 
\item
$\lSdim\T=\displaystyle\lim_{\lambda\to+0}\frac{\log V_{\sigma,\lambda}(\T)}{\log\lambda}=\lim_{t\to-\infty}\frac{h_{\sigma,t}(\S)}{t}$. 
\end{enumerate}
\end{lem}
\begin{pf}
We prove (i). 
Let $A^{(n)}$ be the $\sigma$-semistable factor of $\S^nG$ with $\phi_\sigma(A^{(n)})=\phi^+_\sigma(\S^nG)$. 
We set $C':=\min\{||\gamma||~|~\gamma\in\Lambda\backslash\{0\}\}>0$. 
By the support property we have 
\begin{equation*}
m_{\sigma,\log\lambda}(\S^nG)
\ge \lambda^{\phi^+_\sigma(\S^nG)}|Z_\sigma(v(A^{(n)}))|
> \lambda^{\phi^+_\sigma(\S^nG)}\cdot\frac{1}{C}||v(A^{(n)}))||
\ge \lambda^{\phi^+_\sigma(\S^nG)}\frac{C'}{C}, 
\end{equation*}
which implies $\lambda^{\uSdim\T}\le V_{\sigma,\lambda}(\T)$ by Proposition \ref{lim-of-phase} (i). 
Combining the inequality from Lemma \ref{sigma-vol-scale} (i), 
we have 
\begin{equation*}
\lambda^{\uSdim\T}\le V_{\sigma,\lambda}(\T)\le V_{\sigma,1}(\T)\cdot\lambda^{\uSdim\T},
\end{equation*}
which gives the first equality. 
The second equality in (i) is shown by the definition of $\sigma$-volume.  
The statement (ii) follows from the same argument. 
\qed
\end{pf}
\begin{cor}\label{ob3}
Fix $d\in\RR$. 
Then $\uSdim\T=\lSdim\T=d$ if and only if  
\begin{equation*}
V_{\sigma,\lambda}(\T)=V_{\sigma,1}(\T)\cdot\lambda^d
\end{equation*}
for any $\lambda>0$.
\end{cor}
\begin{pf}
It is clear by Lemma \ref{sigma-vol-scale} and Lemma \ref{serre-sigma-vol}. 
\qed
\end{pf}
\section{Serre dimension and global dimension}\label{sec-stab-Gepner}
In this section, we prove Theorem \ref{Sdim-le-inf} and Theorem \ref{fCY-iff}. 
\begin{lem}\label{key}
For a nonzero object $E\in\T$, we have
\begin{enumerate}
\item
$\phi^+_\sigma(\S E)-\phi^+_\sigma(E)\le\gldim\sigma$ for all $\sigma\in {\rm Stab}_\Lambda(\T)$. 
\item
$\phi^-_\sigma(\S E)-\phi^-_\sigma(E)\le\gldim\sigma$ for all $\sigma\in {\rm Stab}_\Lambda(\T)$. 
\end{enumerate}
\end{lem}
\begin{pf}
We shall prove (i). 
Fix any $\sigma\in {\rm Stab}_\Lambda(\T)$. 
Let $\{A_i\}_{i=1,\cdots,p}$ be the $\sigma$-semistable factors of $E$ with $\phi_\sigma(A_{i-1})>\phi_\sigma(A_i)$, and $(\S E)^+$ the $\sigma$-semistable factor of $\S E$ with $\phi_\sigma((\S E)^+)=\phi^+_\sigma(\S E)$. 
Then, by 
\begin{equation*}
\Hom(E,(\S E)^+)=\Hom((\S E)^+,\S E)^*\neq0, 
\end{equation*}
there exists $i$ such that $\Hom(A_i,(\S E)^+)\neq0$. 
Therefore we have
\begin{equation*}
\phi^+_\sigma(\S E)-\phi^+_\sigma(E)
\le \phi_\sigma((\S E)^+)-\phi_\sigma(A_i)
\le \gldim\sigma. 
\end{equation*}
The proof of (ii) is same by $\Hom(E^-,\S E)=\Hom(E,E^-)^*\neq0$, where $E^-$ is the $\sigma$-semistable factor of $E$ with $\phi_\sigma(E^-)=\phi^-_\sigma(E)$. 
\qed
\end{pf}
The following is a fundamental inequality between the upper Serre dimension and the global dimension.
\begin{thm}\label{Sdim-le-inf}
We have 
$\uSdim\T\le\displaystyle\inf_{\sigma\in {\rm Stab}_\Lambda(\T)}{\gldim\sigma}$. 
\end{thm}
\begin{pf}
Fix $\sigma\in {\rm Stab}_\Lambda(\T)$. 
Since $\Hom(G,\S^nG[-m^-_n])\neq0$ ($n\in\ZZ_{>0}$), we have 
\begin{equation*}
\phi^+_\sigma(\S^nG[-m^-_n])-\phi^-_\sigma(G)\ge0. 
\end{equation*}
It follows from Lemma \ref{key}(i) that
\begin{eqnarray*}
m^-_n
&\le&\phi^+_\sigma(\S^nG)-\phi^-_\sigma(G)\\
&=&(\phi^+_\sigma(\S^nG)-\phi^+_\sigma(\S^{n-1}G))+\cdots+(\phi^+_\sigma(\S G)-\phi^+_\sigma(G))+\phi^+_\sigma(G)-\phi^-_\sigma(G)\\
&\le&n\cdot \gldim\sigma+\phi^+_\sigma(G)-\phi^-_\sigma(G). 
\end{eqnarray*}
Hence we have $\uSdim\T\le\gldim\sigma$. 
\qed
\end{pf}
We give a sufficient condition for the equality between the upper Serre dimension and the infimum of the global dimensions. 
\begin{prop}\label{Sdim-ge-inf}
Fix $s\in\RR$ and assume that, for any $\epsilon\in(0,1)$, there exists $\sigma_{\epsilon}\in {\rm Stab}_\Lambda(\T)$ such that $\S(\P_{\sigma_{\epsilon}}(\phi))\subset\P_{\sigma_{\epsilon}}([\phi+s-\epsilon, \phi+s+\epsilon])$ for all $\phi\in\RR$. 
Then we have
\begin{equation*}
\uSdim\T=\displaystyle\inf_{\sigma\in {\rm Stab}_\Lambda(\T)}{\gldim\sigma}=s. 
\end{equation*}
\end{prop}
\begin{pf}
Fix $\epsilon\in(0,1)$. Then there exists $\sigma_{\epsilon}\in {\rm Stab}_\Lambda(\T)$ such that $\S(\P_{\sigma_{\epsilon}}(\phi))\subset\P_{\sigma_{\epsilon}}([\phi+s-\epsilon, \phi+s+\epsilon])$ for all $\phi\in\RR$ by the assumption. 
For any $\sigma_\epsilon$-semistable object $E\in\T$, we have $\phi^+_{\sigma_\epsilon}(\S E)\le\phi_{\sigma_\epsilon}(E)+s+\epsilon$. 
If $\Hom(E_1,E_2)\neq0$ for $\sigma _\epsilon$-semistable objects $E_1,E_2\in\T$, then $\Hom(E_2,\S E_1)\neq0$. 
This implies
\begin{equation*}
\phi_{\sigma _\epsilon}(E_2)-\phi_{\sigma _\epsilon}(E_1)
\le\phi^+_{\sigma _\epsilon}(\S E_1)-\phi_{\sigma _\epsilon}(E_1)
\le s+\epsilon, 
\end{equation*}
which gives $\gldim\sigma_\epsilon\le s+\epsilon$. 
Hence we have $\displaystyle\inf_{\sigma\in {\rm Stab}_\Lambda(\T)}{\gldim\sigma}\le s$. 

Taking the $\sigma$-semistable factors of a split-generator of $\T$, we get $\sigma_\epsilon$-semistable objects $G_1,\cdots, G_k$ such that $G_1\oplus\cdots\oplus G_k$ is a split-generator of $\T$. 
For each $n\in\ZZ_{>0}$, there exist $i_n, j_n\in\{1,\cdots,k\}$ satisfying $\Hom(G_{i_n}, \S^nG_{j_n}[-m^-_n])\neq0$ ($n\in\ZZ_{>0}$). 
It follows from 
\begin{equation*}
\phi^-_{\sigma_\epsilon}(\S^nG_{j_n})\ge\phi_{\sigma_\epsilon}(G_{j_n})+n(s-\epsilon)
\end{equation*}
that
\begin{eqnarray*}
\gldim\sigma_\epsilon
&\ge&\phi^-_{\sigma_\epsilon}(\S^nG_{j_n}[-m^-_n])-\phi_{\sigma_\epsilon}(G_{i_n})\\
&\ge&n(s-\epsilon)-m^-_n+\phi_{\sigma_\epsilon}(G_{j_n})-\phi_{\sigma_\epsilon}(G_{i_n}). 
\end{eqnarray*}
Hence we have $s-\uSdim\T\le\epsilon$, which implies $s\le\uSdim\T$. 

Therefore the claim follows from Theorem \ref{Sdim-le-inf}. 
\qed
\end{pf}
\begin{cor}\label{Gepner-to-=}
Let $s$ be a real number. 
If $\T$ admits a stability condition $\sigma$ satisfying 
\begin{equation*}
\S(\P_{\sigma}(\phi))=\P_{\sigma}(\phi+s)
\end{equation*}
for all $\phi\in\RR$, then we have
\begin{equation*}
s=\uSdim\T=\displaystyle\inf_{\sigma'\in {\rm Stab}_\Lambda(\T)}{\gldim\sigma'}=\gldim\sigma\ge0. 
\end{equation*}
\end{cor}
\begin{pf}
This immediately follows from Proposition \ref{Sdim-ge-inf}. 
\qed
\end{pf}
The fractional Calabi--Yau condition implies the converse of Corollary \ref{Gepner-to-=}.
\begin{prop}\label{=-to-Gepner}
Suppose that $\T$ is a fractional Calabi--Yau category. 
If there exists $\sigma\in {\rm Stab}_\Lambda(\T)$ such that $\gldim\sigma=\uSdim\T$, 
then $\sigma$ is of Gepner type with respect to $(\S, \uSdim\T)$. 
\end{prop}
\begin{pf}
The Serre functor $\S$ of $\T$ satisfies $\S^{k_1}=[k_2]$ for some $k_1,k_2\in\ZZ$ (in fact $\ZZ_{>0}$). 
Let $E\in\T$ be a $\sigma$-semistable object. 
It follows from Lemma \ref{key}(i) and $\gldim\sigma=\frac{k_2}{k_1}$ that
\begin{eqnarray*}
k_2
&=&\phi_{\sigma}(\S^{k_1}E)-\phi_{\sigma}(E)\\
&=&(\phi_{\sigma}(\S^{k_1}E)-\phi^+_{\sigma}(\S^{k_1-1}E))+(\phi^+_{\sigma}(\S^{k_1-1}E)-\phi^+_{\sigma}(\S^{k_1-2}E))\\
& &+\cdots+(\phi^+_{\sigma}(\S E)-\phi_{\sigma}(E))\\
&\le&k_1\cdot\gldim\sigma=k_2, 
\end{eqnarray*}
which implies $\phi^+_{\sigma}(\S^l E)-\phi^+_{\sigma}(\S^{l-1} E)=\gldim\sigma=\frac{k_2}{k_1}$ for $l=1,\cdots,k_1$. 
By Lemma \ref{key}(ii), we can get $\phi^-_{\sigma}(\S^l E)-\phi^-_{\sigma}(\S^{l-1} E)=\gldim\sigma=\frac{k_2}{k_1}$ for $l=1,\cdots,k_1$ in the same way. 
By induction on $l$, we have $\phi^+_{\sigma}(\S^l E)=\phi^-_{\sigma}(\S^l E)$ for $l=1,\cdots,k_1$, 
that is, $\S^l E$(for all $l\in\ZZ$) is $\sigma$-semistable with $\phi_{\sigma}(\S^l E)=\phi_{\sigma}(E)+l\cdot\frac{k_2}{k_1}$. 

The element $\alpha(\S)\in\Aut_\ZZ(\Lambda)$ is of finite order, which implies 
\begin{equation*}
Z_{\S.\sigma}(v(E))=\exp(-i\pi(\uSdim\T))\cdot Z_{\sigma}(v(E))=Z_{\sigma.(\uSdim\T)}(v(E))
\end{equation*}
for all $\S.\sigma$-semistable objects $E\in\T$. 
Since semistable objects generate $\Lambda$, we have $Z_{\S.\sigma}=Z_{\sigma.(\uSdim\T)}$. 
\qed
\end{pf}
The following is the second main result in this section.
\begin{thm}\label{fCY-iff}
Suppose that $\T$ is a fractional Calabi--Yau category with a stability condition $\sigma$. 
Then, $\gldim\sigma=\uSdim\T$ if and only if $\sigma$ is of Gepner type with respect to $(\S,s)$ for some $s\in\RR$. 
\end{thm}
\begin{pf}
This follows from Corollary \ref{Gepner-to-=} and Proposition \ref{=-to-Gepner}. 
\qed
\end{pf}
\section{Low dimensional triangulated categories}
In this section, we study triangulated categories $\T$ of $\uSdim\T\le1$. 
In Theorem \ref{Sdim<1}, we will classify triangulated categories $\T$ of $\uSdim\T<1$ with a Gepner stability condition. 
\begin{lem}[{\cite[Lemma 3.3]{Qiu}}]\label{ind-semistable}
Let $\sigma$ be a stability condition $\sigma$ on $\T$ with $\gldim\sigma\le1$. 
Then all indecomposable objects in $\T$ are $\sigma$-semistable. 
\end{lem}
\begin{prop}\label{semiconti-inf}
Suppose that $\T$ admits a stability condition $\sigma$ satisfying $\gldim\sigma\le1$. 
Then, for any nonzero admissible triangulated subcategory $\T'$ of $\T$, 
we have
\begin{equation*}
\displaystyle\inf_{\sigma'\in {\rm Stab}_{v(K(\T'))}(\T')}{\gldim\sigma'}\le\inf_{\sigma\in {\rm Stab}_\Lambda(\T)}{\gldim\sigma}. 
\end{equation*} 
\end{prop}
\begin{pf}
Lemma \ref{ind-semistable} implies that all indecomposable objects in $\T$ are $\sigma$-semistable. 
Therefore $\sigma':=(Z_\sigma|_{v(K(\T'))},\{\P_\sigma(\phi)\cap\T'\}_{\phi\in\RR})$ is a stability condition on $\T'$, and clearly satisfies $\gldim\sigma'\le\gldim\sigma$. 
\qed
\end{pf}
\begin{cor}[semicontinuity]\label{semiconti}
Suppose that $\T$ admits a stability condition $\sigma$ satisfying $\uSdim\T=\gldim\sigma\le1$. 
Then the semicontinuity of the Serre dimension holds: for any nonzero admissible triangulated subcategory $\T'$ of $\T$, we have
\begin{equation*}
\uSdim\T'\le\uSdim\T. 
\end{equation*} 
\end{cor}
\begin{pf}
The statement immediately follows from Proposition \ref{semiconti-inf} and Theorem \ref{Sdim-le-inf}. 
\qed
\end{pf}
\subsection{Case of $\uSdim\T<1$}
In this subsection, we prove Theorem \ref{Sdim<1} by studying properties of a stability condition $\sigma$ with $\gldim\sigma < 1$. 
\begin{lem}\label{hereditary}
Let $\sigma$ be a stability condition $\sigma$ on $\T$ with $\gldim \sigma < 1$. 
Then the following hold.
\begin{enumerate}
\item The heart $\P_\sigma((0,1])$ of a bounded t-structure is hereditary.
\item For $\sigma$-semistable objects $E_1, E_2  \in \T$, if $\Hom(E_1,E_2) \neq 0$, we have $\Ext^k(E_1,E_2)=0$ for any nonzero integer $k$.
\item For $E_1 \in \P_\sigma(\phi_1)$ and $E_2 \in \P_\sigma(\phi_2)$, if $\phi_1 \leq \phi_2$, we have $\Ext^k(E_1,E_2)=0$ for any positive integer $k$. 
\end{enumerate}
\end{lem}
\begin{pf}
These immediately follow from the assumption $\gldim \sigma<1$. 
\qed
\end{pf}
\begin{lem}[see also {\cite[Lemma 3.3]{Qiu}}]\label{ind-exc-stab}
Let $\sigma$ be a stability condition $\sigma$ on $\T$ with $\gldim\sigma<1$. For an object $E\in\T$, the following are equivalent.
\begin{enumerate}
\item $E$ is indecomposable.
\item $E$ is exceptional.
\item $E$ is $\sigma$-stable.
\end{enumerate}
\end{lem}
\begin{pf}
The statement (ii)$\Rightarrow$(i) is evident, and (iii)$\Rightarrow$(ii) is clear by $\gldim\sigma<1$. 
We show (i)$\Rightarrow$(iii). 
An object $E$ is $\sigma$-semistable by Lemma \ref{ind-semistable}. 
Assume that $E$ is not $\sigma$-stable. 
Then there is an exact sequence 
\[0 \to E' \to E \to E'' \to 0\]
in $\P_\sigma(\phi_\sigma(E))$ such that $E'$ is nonzero and not isomorphic to $E$. 
By Lemma \ref{hereditary} (iii), we have $\Ext^{1}(E'',E')=0$, hence $E$ is isomorphic to $E'\oplus E''$. 
\qed
\end{pf}
When $\T\simeq\perf(A)$, 
we define $\HH_*(\T):=\HH_*(A)$, where $\HH_*(A)$ is the total space of the Hochschild homology of $A$. 
Since $A$ is smooth proper, we have $\dim_\CC\HH_*(\T)<+\infty$. For a semiorthogonal decomposition $\T=\langle\A,\B\rangle$, $\HH_*(\T)\simeq\HH_*(\A)\oplus\HH(\B)$ holds (see \cite[2.2.8]{Tab}). 

The property of $\gldim\sigma<1$ implies the ``discreteness'' of phases.
\begin{lem}\label{phase-finite}
Let $\sigma$ be a stability condition $\sigma$ on $\T$ with $\gldim \sigma < 1$. For a subset $I \subset \RR$, we put $S_\sigma(I):=\{\phi \in I \mid \P_\sigma(\phi)\neq 0\}$. Then $S_\sigma((0,1])$  is a finite set. 
Moreover, $S_\sigma((0,n))$ is also finite for all positive integers $n$. 
\end{lem}
\begin{pf}
Assume that the set $S_\sigma((0,1])$ is an infinite set.
We can take a monotone increasing sequence $\{\phi_i\}_{i\in \ZZ_{>0}}$ or a monotone decreasing sequence $\{\phi_i\}_{i\in \ZZ_{<0}}$ in $S_\sigma((0,1])$ such that   $\phi_j-\phi_i<1-\gldim\sigma$ for any $i<j$. Take a sequence of $\sigma$-stable objects $\{E_i \in \P_\sigma(\phi_i)\}_i$. We show that $(E_i, E_j)$ is an exceptional pair for $i<j$. By Lemma \ref{ind-exc-stab}, $E_i$ is an exceptional object for any $i$.  Since $\Ext^{\le 0}(E_j, E_i)=0$ for $i<j$, it is enough to show that $\Ext^k(E_j, E_i)=0$ for any positive integer $k$. In fact, it is deduced from
\begin{eqnarray*}
\phi_\sigma(E_i[k])-\phi_\sigma(E_j)
&=&k-(\phi_j-\phi_i)\\
&>&1-(1-\gldim\sigma)\\
&=&\gldim\sigma. 
\end{eqnarray*}
Hence, we have an exceptional collection of infinite length. Since $\dim_\CC\HH_*(\T)<+\infty$, this is a contradiction.
The second statement is deduced from the property $\P_\sigma(\phi)[n]=\P_\sigma(\phi+n)$ for $n\in \ZZ$. \qed
\end{pf}
For a $\CC$-linear category $\C$, denote the set of isomorphism classes of indecomposable objects in $\C$ by $\mathrm{Ind}(\C)$. 
\begin{cor}\label{ind-finite}
Let $\sigma$ be a stability condition $\sigma$ on $\T$ with $\gldim \sigma < 1$. Then $\mathrm{Ind}(\P_\sigma((0,1]))$ is a finite set. 
\end{cor}
\begin{pf}
For any $\phi\in S_\sigma((0,1])$, all objects in $\mathrm{Ind}(\P_\sigma(\phi))$ form a mutually orthogonal exceptional collection by Lemma \ref{hereditary} (iii) and Lemma \ref{ind-exc-stab}. 
Hence $\dim_\CC\HH_*(\T)<+\infty$ implies the finiteness of $\mathrm{Ind}(\P_\sigma(\phi))$. 
By Lemma \ref{ind-semistable}, we have $\mathrm{Ind}(\P_\sigma((0,1]))=\bigcup_{\phi\in S_\sigma((0,1])}\mathrm{Ind}(\P_\sigma(\phi))$. 
The finiteness of $\mathrm{Ind}(\P_\sigma((0,1]))$ follows from Lemma \ref{phase-finite}. 
\qed
\end{pf}
\begin{defn}
A triangulated category $\T$ is {\it connected} if $\T$ has no nontrivial orthogonal decompositions. 
\end{defn}
The property of $\gldim\sigma<1$ also implies the existence of a full strong exceptional collection.
\begin{prop}\label{fec}
Let $\sigma$ be a stability condition $\sigma$ on $\T$ with $\gldim \sigma < 1$. 
Suppose that $\T$ is connected. 
Then $\T$ has a full strong exceptional collection. 
\end{prop}
\begin{pf}
We note that the set $S_\sigma((0,n))$ is finite for all positive integers $n$ by Lemma \ref{phase-finite}. 
Let $\phi_1$ be the minimum number in the set $S_\sigma((0,1])$. Take a $\sigma$-stable object $E_1 \in \P_\sigma(\phi_1)$. 
Assume that $\D_1:=$$^\perp{E}_1 \neq 0$. 
Since $\T$ is connected, there exists an indecomposable object $C\in\D_1$ such that $\Hom(E_1, C)\neq 0$. 
Note that $C$ is indecomposable in $\D$ by the definition of $\D_1$. 
By Lemma \ref{ind-exc-stab} and Definition \ref{stab-def} (iii), $C$ is $\sigma$-stable with $\phi_\sigma(C)>0$. 
Since $\gldim\sigma<1$, the phase $\phi_\sigma(C)$ is lower than $2$, hence $\phi_\sigma(C)\in S_\sigma((0,2))$.  
Define
 \[
 \phi_2:=\mathrm{min}\{\phi_\sigma(C) \in \RR \mid C \in \D_1, \text{\rm $C$ is $\sigma$-stable},  \Hom(E_1, C)\neq 0 \}\in S_\sigma((0,2)).
 \] 
Then $\phi_1 \le \phi_2$ holds.
By the definition of $\phi_2$, 
we can take a $\sigma$-stable object $E_2\in \P_\sigma(\phi_2)$ such that $E_2\in \D_1$ and $\Hom(E_1, E_2)\neq 0$.
By Lemma \ref{hereditary} (ii), we have $\Ext^k(E_1,E_2)=0$ for a nonzero integer $k$, 
hence $(E_1,E_2)$ is a strong exceptional pair by Lemma \ref{ind-exc-stab}. 
Assume that $\D_2:=$$^\perp \langle E_1, E_2 \rangle\neq0$. 
Since $\T$ is connected, there exists an indecomposable object $C\in\D_2$ such that $\Hom(E_1\oplus E_2, C)\neq 0$. 
Note that $C$ is indecomposable in $\D$ by the definition of $\D_2$. 
By Lemma \ref{ind-exc-stab} and Definition \ref{stab-def} (iii), $C$ is $\sigma$-stable with $\phi_\sigma(C)>0$. 
Since $\gldim\sigma<1$, the phase $\phi_\sigma(C)$ is lower than $3$, hence $\phi_\sigma(C)\in S_\sigma((0,3))$.
Define
\[ \phi_3:=\mathrm{min}\{\phi_\sigma(C) \in \RR \mid C \in\D_2, \text{\rm $C$ is $\sigma$-stable},  \Hom(E_1\oplus E_2, C)\neq 0 \}\in S_\sigma((0,3)). 
\]
Then $\phi_1 \le \phi_2 \le \phi_3$ holds. 
Take a $\sigma$-stable object $E_3\in \P_\sigma(\phi_3)$ such that $E_3\in \D_2$ and $\Hom(E_1\oplus E_2, E_3)\neq 0$. 
By Lemma \ref{hereditary} (iii), we have $\Ext^{>0}(E_i,E_3)=0$ for $i=1,2$. 
By the minimality of $\phi_3$, we have $\Ext^{<0}(E_i,E_3)=0$ for $i=1,2$. 
Therefore $(E_1,E_2,E_3)$ is a strong exceptional collection by Lemma \ref{ind-exc-stab}. 

We continue this procedure until $\D_n:=$$^\perp\langle E_1, E_2, \cdot \cdot \cdot, E_n\rangle$ becomes zero. 
Due to $\dim_\CC\HH_*(\T) < \infty$, such positive integer $n$ exists.
\qed
\end{pf}
We recall the notion of locally finiteness of triangulated categories.
\begin{defn}
A triangulated category $\T$ is {\it locally finite} if for any object $F \in \T$ there are only finitely many isomorphism classes of indecomposable objects $E$ such that $\Hom_\T(E,F)\neq0$. 
\end{defn}
\begin{prop}[Auslander, Happel and Beligiannis (see {\cite[Prop.2.3, Examples(2)]{Kra}})]\label{dynkin}
Let $B$ be a finite dimensional $\CC$-algebra such that $\D^b(B)$ is connected. 
Then $\D^b(B)$ is locally finite if and only if there is an equivalence $\D^b(B) \simeq \D^b(Q)$ for some Dynkin quiver $Q$. 
\end{prop}
The following is the main theorem in this section. 
\begin{thm}\label{Sdim<1}
Suppose that $\T$ is connected. 
The following are equivalent.
\begin{enumerate}
\item
$\T$ is equivalent to $\D^b(Q)$ for some Dynkin quiver $Q$. 
\item
$\uSdim\T<1$ and there exists a Gepner type stability condition on $\T$. 
\item
There exists $\sigma\in {\rm Stab}_\Lambda(\T)$ with $\gldim\sigma<1$. 
\end{enumerate}
\end{thm}
\begin{pf}
We shall consider (i)$\Rightarrow$(ii). By Example \ref{example-dimension}(iii), (iv), we have $\uSdim\T<1$. 
By Example \ref{KST-Gepner}, $\T$ admits a Gepner type stability condition $\sigma_G$. 
(ii)$\Rightarrow$(iii) follows from Corollary \ref{Gepner-to-=}. 
Then we shall show (iii)$\Rightarrow$(i). By Proposition \ref{fec}, we have a full strong exceptional collection $\T=\langle E_1,\cdots, E_n \rangle$. 
Then we have $\T \simeq D^b(B)$, where $B:=\mathrm{End}_\T(\bigoplus_{i=1}^{n}E_i)$. 
Since $\T$ is locally finite by Corollary \ref{ind-finite} and $\gldim\sigma<1$, 
Proposition \ref{dynkin} gives an equivalence $\T\simeq\D^b(Q)$ for some Dynkin quiver $Q$.
\qed
\end{pf}
\subsection{Case of $\uSdim\T=1$}
We calculate the global dimensions in the case of curves in this subsection. 

For an acyclic quiver $Q$, not of Dynkin type, $\uSdim\D^b(Q)(=1)$ is equal to the global dimension of some stability condition by the following theorem due to Qiu.
\begin{thm}[{\cite[Theorem 5.2]{Qiu}}]\label{gldimQ}
Suppose that $Q$ is an acyclic quiver, not of Dynkin type. 
Then there exists a stability condition $\sigma$ on $\D^b(Q)\in {\rm Stab}_K(\D^b(Q))$, such that
\begin{equation*}
\displaystyle\inf_{\sigma'\in {\rm Stab}_K(\D^b(Q))}{\gldim\sigma'}=\gldim\sigma
=1. 
\end{equation*}
\end{thm}
Next, we treat derived categories of smooth projective curves. 

\begin{thm}[{\cite[Theorem 2.7]{Mac}}]\label{stab-curve}
Let $C$ be a smooth projective curve of genus $g \geq 1$. 
For $\beta\in\RR$ and $H\in\RR_{>0}$, we define a group homomorphism $Z_{\beta,H}: \N(\D^b(C))\to\CC$ by
\begin{equation*}
Z_{\beta,H}(E):=-\mathrm{deg}(E)+(\beta+iH)\cdot \mathrm{rk}(E).
\end{equation*}
The pair $\sigma_{\beta,H}:=(Z_{\beta,H}, {\rm Coh}(C))$ is a numerical stability condition on $\D^b(C)$.
Moreover, the map 
\[ \mathbb{H}\times\CC \to \mathrm{Stab}_\N(\D^b(C)),(\beta+iH, \mu) \mapsto \sigma_{\beta,H}.\mu \]
is an isomorphism, where $\mathbb{H}$ is the upper half plane.
\end{thm}

We prepare the following lemma for applying Proposition \ref{Sdim-ge-inf}.
\begin{lem}\label{largevolume}
Let $C$ be a smooth projective curve of genus $g\geq2$. 
Then for any $\epsilon\in(0,1)$, there exists $\sigma_\epsilon\in {\rm Stab}_\N(\D^b(C))$ such that, for all $\phi\in\RR$, 
\begin{equation*}
\S(\P_{\sigma_{\epsilon}}(\phi))\subset\P_{\sigma_{\epsilon}}([\phi+1-\epsilon, \phi+1+\epsilon]). 
\end{equation*}
\end{lem}
\begin{pf}
Recall that $\S=-\otimes \omega_C[1]$.
Fix $\epsilon\in(0,1)$. 
Since the function $\mathrm{arccot}:\RR\to(0,\pi)$ is uniformly continuous, there exists $H>0$ such that 
\begin{equation*}
\frac{1}{\pi}\left|\mathrm{arccot}\Bigl(x-\frac{2g-2}{H}\Bigr)-\mathrm{arccot}(x)\right|<\epsilon. 
\end{equation*}
for all $x\in\RR$. 
Define $\sigma_\epsilon:=\sigma_{0,H}\in {\rm Stab}_\N(\D^b(C))$. 
Let $E$ be a $\sigma_\epsilon$-semistable object. When $\mathrm{rk}(E)=0$, we have $\phi_{\sigma_\epsilon}(E\otimes \omega_C)-\phi_{\sigma_\epsilon}(E)=0$. 
Then we treat the case of $\mathrm{rk}(E) \neq 0$. 
Note that $Z_{0,H}(E\otimes \omega_C)=Z_{0,H}(E)-(2g-2)\mathrm{rk}(E)$. 
Therefore we have 
\begin{eqnarray*}
|\phi_{\sigma_\epsilon}(E \otimes \omega_C)-\phi_{\sigma_\epsilon}(E)|
&=&\frac{1}{\pi}\left|\mathrm{arg}(Z_{0,H}(E\otimes \omega_C))-\mathrm{arg}(Z_{0,H}(E))\right|\\
&=&\frac{1}{\pi}\left|\mathrm{arg}\Bigl(-\frac{\mathrm{deg}(E)}{H\cdot \mathrm{rk}(E)}-\frac{2g-2}{H}+i\Bigr)-\mathrm{arg}\Bigl(-\frac{\mathrm{deg}(E)}{H\cdot \mathrm{rk}(E)}+i\Bigr)\right|\\
&=&\frac{1}{\pi}\left|\mathrm{arccot}\Bigl(-\frac{\mathrm{deg}(E)}{H\cdot \mathrm{rk}(E)}-\frac{2g-2}{H}\Bigr)-\mathrm{arccot}\Bigl(-\frac{\mathrm{deg}(E)}{H\cdot \mathrm{rk}(E)}\Bigr)\right|\\
&<&\epsilon, 
\end{eqnarray*}
which gives the claim. 
\qed
\end{pf}
The behavior of (the infimum of) the global dimensions in the case of smooth projective curves of genus $\ge2$ is different from the quiver case.
\begin{thm}\label{curve-inf}
Let $C$ be a smooth projective curve of genus $g$.
The following holds. 
\begin{enumerate}
\item
If $g=0$, then there is a stability condition on $\D^b(C)$ such that $\gldim \sigma=1$.
\item
If $g=1$, then we have $\gldim \sigma=1$ for any stability condition $\sigma \in {\rm Stab}_\N(\D^b(C))$.
\item
If $g \geq 2$, then $\gldim \sigma>1$ holds for any stability condition  $\sigma \in {\rm Stab}_\N(\D^b(C))$ and $\displaystyle\inf_{\sigma\in {\rm Stab}_\N(\D^b(C))}{\gldim\sigma}=1$.
\end{enumerate}
\end{thm}
\begin{pf}
If $g=0$, $C$ is the projective line and $\D^b(C)$ is equivalent to $\D^b(K_2)$, where $K_2$ is the Kronecker quiver with two arrows. 
By Theorem \ref{gldimQ}, we have (i). 
If $g=1$, $C$ is an elliptic curve. 
Since $\D^b(C)$ is a one dimensional Calabi--Yau category, we have (ii) by Corollary \ref{Gepner-to-=}. 

We prove the first statement in (iii). 
Take a stability condition $\sigma\in {\rm Stab}_\N(\D^b(C))$. 
By Theorem \ref{stab-curve}, we may assume that $\sigma=\sigma_{\beta, H}$ for $\beta \in \RR$ and  $H \in \RR_{>0}$. 
Note that $\mathcal{O}_C$ and $\omega_C[1]$ are $\sigma_{\beta,H}$-semistable objects.
Since $\Hom_{\D^b(C)}(\mathcal{O}_C, \omega_C[1])=\CC$, we have
\begin{equation*}
\gldim \sigma_{\beta,H} \geq 1+\phi_{\sigma_{\beta,H}}(\omega_C)-\phi_{\sigma_{\beta,H}}(\mathcal{O}_C). 
\end{equation*}

Note that $Z_{\beta,H}(\mathcal{O}_C)=\beta+iH$ and $Z_{\beta,H}(\omega_C)=\beta-(2g-2)+iH$. 
By $g\geq2$, we have $\phi_{\sigma_{\beta,H}}(\omega_C)>\phi_{\sigma_{\beta,H}}(\mathcal{O}_C)$, which implies $\gldim \sigma>1$. 
The second statement in (iii) is clear by Lemma \ref{largevolume} and Proposition \ref{Sdim-ge-inf}. 
\qed
\end{pf}
Theorem \ref{Sdim<1} completely classifies $\T$ of $\uSdim\T<1$ with a Gepner type stability condition. 
Therefore the next step is the classification of the case of $\uSdim\T=1$, and Corollary \ref{semiconti} might be useful. 
We hope to return to this topic in future research. 

\end{document}